\documentclass{elsart}

\usepackage{amssymb}
\usepackage{amsmath}
\usepackage[english]{babel}

\newtheorem{theorem}{Theorem}

\newtheorem{corollary}[theorem]{Corollary}

\newtheorem{lemma}[theorem]{Lemma}

\newtheorem{remark}[theorem]{Remark}
\begin{document}
\begin{frontmatter}

\title{Contributions of non-extremum critical points to the semi-classical trace formula.}

\author{Brice Camus}
\thanks[label2]{This work was supported by the
\textit{IHP-Network}, \textit{Analysis\&Quantum} reference
HPRN-CT-2002-00277. We thank Yves Colin de Verdi\`ere that
indicated us the result of Yorke.}
\address{Mathematisches Institut der Ludwig-Maximilians-Universit\"{a}t M\"{u}nchen.\\
Theresienstra\ss{}e 39, 80333 Munich, Germany.\newline%
Email : camus@mathematik.uni-muenchen.de}

\begin{abstract}
We study the semi-classical trace formula at a critical energy
level for a $h$-pseudo-differential operator on $\mathbb{R}^{n}$
whose principal symbol has a totally degenerate critical point for
that energy. We compute the contribution to the trace formula of
isolated non-extremum critical points under a condition of "real
principal type". The new contribution to the trace formula is
valid for all time in a compact subset of $\mathbb{R}$ but the
result is modest since we have restrictions on the dimension.
\end{abstract}

\end{frontmatter}

\section{Introduction.}
The semi-classical trace formula for a self-adjoint
$h$-pseudo-differential operator $P_{h}$, or more generally
$h$-admissible (see \cite{[Rob]}), studies the asymptotic
behavior, as $h$ tends to 0, of the spectral function :
\begin{equation}
\gamma (E,h,\varphi )=\sum\limits_{|\lambda _{j}(h)-E|\leq
\varepsilon }\varphi (\frac{\lambda _{j}(h)-E}{h}),  \label{Def
trace}
\end{equation}
where the $\lambda _{j}(h)$ are the eigenvalues of $P_{h}$. Here
we suppose that the spectrum is discrete in $[E-\varepsilon
,E+\varepsilon ]$, some sufficient conditions for this will be
given below. If $p_0$ is the principal symbol of $P_{h}$ we recall
that an energy $E$ is regular when $\nabla p_0(x,\xi )\neq 0$ on
the energy surface $\Sigma _{E}=\{(x,\xi )\in
T^{\star}\mathbb{R}^{n}\text{ }/\text{ }p_0(x,\xi )=E\}$ and
critical when it is not regular.

It is well known that the asymptotics of (\ref{Def trace}), as $h$
tends to 0, is closely related to the closed trajectories of the
Hamiltonian flow of $p$ on the surface $\Sigma_{E}$, i.e. :
\begin{equation*}
\lim_{h\rightarrow 0}\gamma (E,h,\varphi )\rightleftharpoons
\{(t,x,\xi)\in\mathbb{R}\times \Sigma _{E} \text{ / }
\Phi_{t}(x,\xi)=(x,\xi)\},
\end{equation*}
where :
\begin{equation*}
\Phi_{t}=\mathrm{exp}(tH_{p_0}):T^{*}\mathbb{R}^{n}\rightarrow
T^{*}\mathbb{R}^{n},
\end{equation*}
and $H_{p_0}=\partial _{\xi }p_0.\partial_{x} -\partial
_{x}p_0.\partial_{\xi}$ is the Hamiltonian vector field attached
to $p_0$.

When $E$ is a regular energy a non-exhaustive list of references
concerning this subject is Gutzwiller \cite{GUT}, Balian and Bloch
\cite{BB} for the physic literature  and for a mathematical point
of view Helffer and Robert \cite{HR}, Brummelhuis and Uribe
\cite{BU}, Paul and Uribe \cite{PU}, and more recently Combescure
et al. \cite{CRR}, Petkov and Popov \cite{P-P}, Charbonnel and
Popov \cite{C-P}.

When one drops the assumption that $E$ be a regular value, the
behavior of (\ref{Def trace}) will depend on the nature of the
singularities of $p$ on $\Sigma_{E}$ which can be complicated. The
semi-classical trace formula for a non-degenerate critical energy,
that is such that the critical-set $\mathbb{\frak{C}}(p_0)
=\{(x,\xi )\in T^{\ast }\mathbb{R}^{n}\text{ }/\text{ }dp_0(x,\xi
)=0\}$ is a compact $C^{\infty }$ manifold with a Hessian
$d^{2}p_0$ transversely non-degenerate along this manifold has
been studied first by Brummelhuis et al. in \cite{BPU}. They
studied this question for quite general operators but for some
''small times'', that is they have supposed that the support of
$\hat{\varphi}$ is contained in such a small neighborhood of the
origin that the only period of the linearized flow in
$\rm{supp}(\hat{\varphi})$ is 0. Later, Khuat-Duy in
\cite{KhD,KhD1} has obtained the contributions of the non-zero
periods of the linearized flow with the assumption that
$\rm{supp}(\hat{\varphi})$ is compact, but for Schr\"{o}dinger
operators with symbol $\xi ^{2}+V(x)$ and a non-degenerate
potential $V$. Our contribution to this subject was to compute the
contributions of the non-zero periods of the linearized flow for
some more general operators, always with $\hat{\varphi}$ of
compact support and under some geometrical assumptions on the flow
(see \cite{Cam0} or \cite{Cam}). Finally, in \cite{Cam1} we have
obtained the contributions to the semi-classical trace formula of
totally degenerate extremum and the objective of this article is
to obtain a generalization when one drop the extremum condition.

Basically, the asymptotics of (\ref{Def trace}) can be expressed
in terms of oscillatory integrals whose phases are related to the
flow of $p_0$ on $\Sigma_E$. When $(x_0,\xi_0)$ is a critical
point of $p_0$, it is well known that the relation
$\mathrm{Ker}(d_{x,\xi}\Phi_{t}(x_0,\xi_0)-\mathrm{Id})\neq \{
0\}$ leads to the study of degenerate oscillatory integrals. Here
we examine the case of a totally degenerate energy, that is such
that the Hessian matrix at our critical point is zero. Hence, the
linearized flow for such a critical point satisfies $d_{x,\xi
}\Phi _{t}(x_{0},\xi _{0})=\mathrm{Id}$, for all $t\in\mathbb{R}$
and the oscillatory integrals we have to consider are totally
degenerate.

The core of the proof lies in establishing suitable normal forms
for our phase functions and in a generalization of the stationary
phase formula for these normal forms. The construction proposed
for the normal forms is independent of the dimension but the
asymptotic expansion of the related oscillatory integrals depends
on the dimension and on the order of the singularity at the
critical point. This explains why the main result is stated with a
restriction on the dimension.

\section{Hypotheses and main result.}
Let $P_{h}=Op_{h}^{w}(p(x,\xi ,h))$ be a $h$-pseudodifferential
operator, obtained by Weyl quantization, in the class of
$h$-admissible operators with symbol $p(x,\xi ,h) \sim \sum
h^{j}p_{j}(x,\xi )$. This means that there exists sequences
$p_{j}\in \Sigma _{0}^{m}(T^{\ast }\mathbb{R} ^{n})$ and
$R_{N}(h)$ such that :
\begin{equation*}
P_{h}=\sum\limits_{j<N}h^{j}p_{j}^{w}(x,hD_{x})+h^{N}R_{N}(h),\text{
} \forall N\in \mathbb{N},
\end{equation*}
where $R_{N}(h)$ is a bounded family of operators on
$L^{2}(\mathbb{R}^{n})$, for $h\leq h_{0}$, and :
\begin{equation*}\Sigma _{0}^{m}(T^{\ast }\mathbb{R}^{n})=\{a:T^{\ast }\mathbb{R}%
^{n}\rightarrow \mathbb{C},\text{ }\sup|\partial _{x}^{\alpha
}\partial _{\xi }^{\beta }a(x,\xi )|<C_{\alpha ,\beta }m(x,\xi
),\text{ }\forall \alpha ,\beta \in \mathbb{N}^{n}\},
\end{equation*}
where  $m$ is a tempered weight on $T^{\ast }\mathbb{R}^{n}$. For
a detailed exposition on $h$-admissible operators we refer to the
book of Robert \cite{[Rob]}. In particular, $p_{0}(x,\xi )$ is the
principal symbol of $P_{h}$ and $p_{1}(x,\xi )$ the sub-principal
symbol. Let be $\Phi _{t}=\textrm{exp} (tH_{p_0}):T^{\ast
}\mathbb{R}^{n}\rightarrow T^{\ast }\mathbb{R}^{n}$, the
Hamiltonian flow of $H_{p_0}=\partial _{\xi }p_0
.\partial_{x}-\partial _{x}p_0 .\partial_{\xi }$.

We study here the asymptotics of the spectral function :
\begin{equation}
\gamma (E_{c},h)=\sum\limits_{\lambda _{j}(h)\in \lbrack
E_{c}-\varepsilon ,E_{c}+\varepsilon ]}\varphi (\frac{\lambda
_{j}(h)-E_{c}}{h}), \label{Objet trace}
\end{equation}
under the hypotheses $(H_{1})$ to $(H_{4})$ given below.
\medskip\\
$(H_{1})$\textit{ The symbol of }$P_{h}$ \textit{is real and there
exists } $\varepsilon_{0}>0$ \textit{ such that the set }%
$p_0^{-1}([E_c-\varepsilon_{0},E_c+\varepsilon_{0}])$\textit{\ is
compact in $T^*\mathbb{R}^n$.}
\begin{remark} \rm{By Theorem 3.13 of \cite{[Rob]} the spectrum $\sigma (P_{h})\cap
[E_{c}-\varepsilon ,E_{c}+\varepsilon ]$ is discrete and consists
in a sequence $\lambda _{1}(h)\leq \lambda _{2}(h)\leq ...\leq
\lambda _{j}(h)$ of eigenvalues of finite multiplicities, if
$\varepsilon$ and $h$ are small enough.}
\end{remark}
To simplify notations we write $z=(x,\xi)$ for any point of the
phase space.\medskip\\
$(H_{2})$\textit{ On }$\Sigma
_{E_{C}}=p_0^{-1}(\{E_{c}\})$\textit{, }$p_0$\textit{ has a unique
critical point }$z_{0}=(x_{0},\xi _{0})$\textit{ and near }$z_{0}$
:
\begin{equation*}
p_0(z)=E_{c}+\sum\limits_{j=k}^{N}\mathfrak{p}_{j}(z)+\mathcal{O}(||(z-z_{0})||^{N+1}),\text{
} k>2,
\end{equation*}
\textit{where the functions }$\mathfrak{p}_{j}$\textit{\ are homogeneous of degree }$j$%
\textit{\ in }$z-z_{0}.$\medskip\\
$(H_{3})$\textit{ We have }$\hat{\varphi}\in C_{0}^{\infty
}(\mathbb{R}).$\ \medskip\\
Since we are interested in the contribution of the fixed point
$z_{0}$, to understand the new phenomenon it is suffices to study
:
\begin{equation}
\gamma _{z_{0}}(E_{c},h)=\frac{1}{2\pi }\mathrm{Tr}\int\limits_{\mathbb{R}}e^{i%
\frac{tE_{c}}{h}}\hat{\varphi}(t)\psi ^{w}(x,hD_{x})\mathrm{exp}(-\frac{i}{h}%
tP_{h})\Theta (P_{h})dt.
\end{equation}
Here $\Theta $ is a function of localization near the critical
energy surface $\Sigma_{E_c}$ and $\psi \in C_{0}^{\infty
}(T^{\ast }\mathbb{R}^{n})$ has an appropriate support near
$z_{0}$. Rigorous justifications are given in section 3 for the
introduction of $\Theta (P_{h})$ and in section 4 for $\psi
^{w}(x,hD_{x})$.\medskip\\
In \cite{Cam1} it was proven that :
\begin{theorem}\label{Old theo}
Under hypotheses $(H_{1})$ to $(H_{3})$ and if $z_{0}$ is a local
extremum of the principal symbol $p_0$ we have
\begin{equation*}
\gamma _{z_{0}}(E_{c},h)\sim
h^{\frac{2n}{k}-n}(\sum\limits_{j=0}^{N}\Lambda _{j,k}(\varphi
)h^{\frac{j}{k}}+\mathcal{O}(h^{\frac{N+1}{k}})), \text{ as
}h\rightarrow 0,
\end{equation*}
where the $\Lambda _{j,k}$ are some distributions and the leading
coefficient is given by :
\begin{equation}
\Lambda _{0,k}(\varphi ) =\frac{1}{k} \left\langle\varphi
(t+p_{1}(z_{0})),t_{z_0}^{\frac{2n-k}{k}}\right\rangle
\frac{1}{(2\pi
)^{n}}\int\limits_{\mathbb{S}^{2n-1}}|\mathfrak{p}_{k}(\theta
)|^{-\frac{2n}{k}}d\theta,
\end{equation}
with $t_{z_0}=\mathrm{max}(t,0)$ if $z_0$ is a minimum and
$t_{z_0}=\mathrm{max}(-t,0)$ for a maximum.
\end{theorem}
To obtain a generalization of Theorem \ref{Old theo}, when the
critical point $z_0$
is not a local extremum, we consider the classical hypothesis : \medskip\\
$(H_{4})$ \textit{We have }$\nabla \mathfrak{p}_{k}\neq 0$\textit{
on the set }$C_{\mathfrak{p}_k}=\{ \theta \in
\mathbb{S}^{2n-1}\text{ / } \mathfrak{p}_{k}(\theta )=0
\}$\textit{.}
\begin{remark}
\rm{We would like to emphasize that, contrary to the case of a
local extremum, the critical point $z_0$ is not necessarily
isolated on $\Sigma_{E_c}$. This imposes to study the classical
dynamic in a micro-local neighborhood of $z_0$. Moreover, by
homogeneity $(H_{4})$ implies that $\nabla \mathfrak{p}_{k}\neq 0$
on the cone $\{(x,\xi)\neq 0\text{ / }\mathfrak{p}_k(x,\xi)=0\}$.
This allows to define the Liouville measure on the set
$C_{\mathfrak{p}_k}$.}
\end{remark}
Then, the new contribution to the trace formula is given by :
\begin{theorem}
\label{Main1} Under hypotheses $(H_{1})$ to $(H_{4})$ and if
$k>2n$ we have :
\begin{equation*}
\gamma _{z_{0}}(E_{c},h)\sim h^{\frac{2n}{k}-n}\Lambda
_{0,k}(\varphi )+\mathcal{O}(h
^{\frac{2n+1}{k}-n}\mathrm{log}^2(h)),\text{ as }h\rightarrow 0,
\end{equation*}
where the leading coefficient is given by :
\begin{gather}
\Lambda _{0,k}(\varphi )
=\frac{1}{k}%
\left ( \left\langle
|t|_{+}^{\frac{2n}{k}-1},\varphi(t+p_1(z_0))\right\rangle
\frac{1}{(2\pi)^{n}}\int\limits_{\mathbb{S}^{2n-1}\cap \{
\mathfrak{p}_k \geq 0\}} |\mathfrak{p}_k
(\theta) |^{-\frac{2n}{k}}d\theta\right. \notag \\
\left. + \left\langle
|t|_{-}^{\frac{2n}{k}-1},\varphi(t+p_1(z_0))\right\rangle
\frac{1}{(2\pi)^{n}}\int\limits_{\mathbb{S}^{2n-1}\cap \{
\mathfrak{p}_k \leq 0\}} |\mathfrak{p}_k (\theta)
|^{-\frac{2n}{k}}d\theta \right ). \label{resultat1}
\end{gather}
If $k=2n$ a similar result is :
\begin{equation*}
\gamma _{z_{0}}(E_{c},h)\sim h^{1-n}\mathrm{log}(h)\Lambda
(\varphi )+\mathcal{O}(h ^{1-n}),\text{ as }h\rightarrow 0,
\end{equation*}
with leading coefficient :
\begin{equation} \Lambda (\varphi )
=\frac{1}{(2\pi
)^n}\hat{\varphi}(0)\mathrm{LVol}(C_{\mathfrak{p}_k}),
\end{equation}
where $\mathrm{LVol}$ is the Liouville measure attached to
$\mathfrak{p}_k$ restricted to the sphere :
\begin{equation*}
\mathrm{LVol}(C_{\mathfrak{p}_k})=\int\limits_{C_{\mathfrak{p}_k}}
dL_{\mathfrak{p}_k}(\theta),
\end{equation*}
with $dL_{\mathfrak{p}_k}(\theta)\wedge
d_{\theta}\mathfrak{p}_k(\theta)=d\theta$ on $C_{\mathfrak{p}_k}$.
\end{theorem}
\begin{remark}
\rm{ In Eq. (\ref{resultat1}) and under the assumption on the
dimension both terms are well defined since for the first term
$|t|^{\frac{2n}{k}-1}$ is locally integrable and for the second
term the integral is convergent.}
\end{remark}
\section{Oscillatory representation.}
Let be $\varphi \in \mathcal{S}(\mathbb{R})$ with
$\hat{\varphi}\in C_{0}^{\infty }(\mathbb{R})$, we recall that :
\begin{equation*} \gamma (E_{c},h)
=\sum\limits_{\lambda _{j}(h)\in I_{\varepsilon }}\varphi
(\frac{\lambda _{j}(h)-E_{c}}{h}), \text{ } I_{\varepsilon }
=[E_{c}-\varepsilon ,E_{c}+\varepsilon ],
\end{equation*}
with $p_{0}^{-1}(I_{\varepsilon_{0}})$ compact in $T^{\ast
}\mathbb{R}^{n}$. By Proposition 3.13 of \cite{[Rob]} the spectrum
of $P_{h}$ is discrete in $I_{\varepsilon }$ for $h>0$ small
enough and $\varepsilon<\varepsilon_{0}$. Now, we localize near
the critical energy $E_{c}$ with a cut-off function $\Theta \in
C_{0}^{\infty }(]E_{c}-\varepsilon ,E_{c}+\varepsilon \lbrack )$,
such that $\Theta =1$ near $E_{c}$ and $0\leq \Theta \leq 1$ on
$\mathbb{R}$. The associated decomposition is :
\begin{equation*}
\gamma (E_{c},h)=\gamma _{1}(E_{c},h)+\gamma _{2}(E_{c},h),
\end{equation*}
with :
\begin{equation}
\gamma _{1}(E_{c},h)=\sum\limits_{\lambda _{j}(h)\in
I_{\varepsilon }}(1-\Theta )(\lambda _{j}(h))\varphi
(\frac{\lambda _{j}(h)-E_{c}}{h}),
\end{equation}
\begin{equation}
\gamma _{2}(E_{c},h)=\sum\limits_{\lambda _{j}(h)\in
I_{\varepsilon }}\Theta (\lambda _{j}(h))\varphi (\frac{\lambda
_{j}(h)-E_{c}}{h}).
\end{equation}
The asymptotic behavior of $\gamma _{1}(E_{c},h)$ is classical and
is given by :
\begin{lemma}
$\gamma _{1}(E_{c},h)=\mathcal{O}(h^{\infty })$, as $h\rightarrow
0$.\label{S1(h)=Tr}
\end{lemma}
For a proof see e.g. \cite{Cam1}.\medskip\\
Consequently, for the study of $\gamma (E_{c},h)$ modulo
$\mathcal{O}(h^{\infty })$, we have only to consider the quantity
$\gamma _{2}(E_{c},h)$. By inversion of the Fourier transform we
obtain the identity :
\begin{equation*}
\Theta (P_{h})\varphi (\frac{P_{h}-E_{c}}{h})=\frac{1}{2\pi}\int\limits_{%
\mathbb{R}}e^{i\frac{tE_{c}}{h}}\hat{\varphi}(t)\mathrm{exp}(-\frac{i}{h}%
tP_{h})\Theta (P_{h})dt.
\end{equation*}
Since the trace of the left hand-side is exactly $\gamma
_{2}(E_{c},h)$, we have :
\begin{equation}
\gamma _{2}(E_{c},h)
=\frac{1}{2\pi}\mathrm{Tr}\int\limits_{\mathbb{R}}e^{i\frac{tE_{c}}{h}}
\hat{\varphi}(t)\mathrm{exp}(-\frac{i}{h}tP_{h})\Theta (P_{h})dt,
\label{Trace S2(h)}
\end{equation}
and with Lemma \ref{S1(h)=Tr} this gives :
\begin{equation*}
\gamma (E_{c},h)=\frac{1}{2\pi }\mathrm{Tr}\int\limits_{\mathbb{R}}e^{i%
\frac{tE_{c}}{h}}\hat{\varphi}(t)\mathrm{exp}(-\frac{i}{h}tP_{h})\Theta
(P_{h})dt+\mathcal{O}(h^{\infty }).
\end{equation*}
Let be $U_{h}(t)=\mathrm{exp}(-\frac{it}{h}P_{h})$, the evolution
operator. For each integer $N$ we can approximate $U_{h}(t)\Theta
(P_{h})$, modulo $\mathcal{O}(h^{N})$, by a Fourier
integral-operator, or FIO, depending on a parameter $h$. Let
$\Lambda$ be the Lagrangian manifold associated to the flow of
$p_0$, i.e. :
\begin{equation*}
\Lambda =\{(t,\tau ,x,\xi ,y,\eta )\in T^{\ast }\mathbb{R}\times T^{\ast }%
\mathbb{R}^{n}\times T^{\ast }\mathbb{R}^{n}:\tau =p(x,\xi
),\text{ }(x,\xi )=\Phi _{t}(y,\eta )\}.
\end{equation*}
\begin{theorem}
The operator $U_{h}(t)\Theta (P_{h})$ is $h$-FIO associated to
$\Lambda$, there exist $U_{\Theta ,h}^{(N)}(t)$ with integral
kernel in $I(\mathbb{R}^{2n+1},\Lambda )$ and $R_{h}^{(N)}(t)$
bounded, with a $L^{2}$-norm uniformly bounded for $0<h\leq 1$ and
$t$ in a compact subset of $\mathbb{R}$, such that
$U_{h}(t)\Theta(P_{h})=U_{\Theta,h}^{(N)}(t)+h^{N}R_{h}^{(N)}(t)$.
\end{theorem}
We refer to Duistermaat \cite{DUI1} for a proof of this theorem.
\begin{remark}
\rm{By a theorem of Helffer and Robert, see e.g. \cite{[Rob]},
Theorem 3.11 and Remark 3.14, $\Theta (P_{h})$ is an
$h$-admissible operator with a symbol supported in
$p_0^{-1}(I_{\varepsilon }).$ This allows us to consider only
oscillatory-integrals with compact support.}
\end{remark}
For the control of the remainder, associated to $R_{h}^{(N)}(t)$,
we use :
\begin{corollary}
Let be $\Theta _{1}\in C_{0}^{\infty }(\mathbb{R})$ such that
$\Theta _{1}=1$ on $\rm{supp}(\Theta )$ and $\rm{supp}(\Theta
_{1})\subset I_{\varepsilon }$, then $\forall N\in \mathbb{N}$ :
\begin{equation*}
\mathrm{Tr}(\Theta (P_{h})\varphi (\frac{P_{h}-E_{c}}{h}))=\frac{1}{2\pi }%
\mathrm{Tr}\int\limits_{\mathbb{R}}\hat{\varphi}(t)e^{\frac{i}{h}%
tE_{c}}U_{\Theta ,h}^{(N)}(t)\Theta
_{1}(P_{h})dt+\mathcal{O}(h^{N}).
\end{equation*}
\end{corollary}
For a proof, see e.g. \cite{Cam1}.\medskip\\
If $(x_{0},\xi _{0})\in \Lambda $ and if $\varphi =$ $\varphi
(x,\theta )\in C^{\infty }(\mathbb{R}^{k}\times \mathbb{R}^{N})$
parameterizes $\Lambda $ in a sufficiently small neighborhood $U$
of $(x_{0},\xi _{0})$ then for each $u_{h}\in
I(\mathbb{R}^{k},\Lambda )$ and $\chi \in C_{0}^{\infty }(T^{\ast
}\mathbb{R}^{k})$, $\rm{supp}(\chi )\subset U,$ there exists a
sequence of amplitudes $a_{j}=a_{j}(x,\theta )\in C_{0}^{\infty }(\mathbb{R}%
^{k}\times \mathbb{R}^{N})$ such that for all $L\in\mathbb{N}$ :
\begin{equation}
\chi ^{w}(x,hD_{x})u_{h}=\sum\limits_{-d\leq j<L}h^{j}I(a_{j}e^{\frac{i}{h}%
\varphi })+\mathcal{O}(h^{L}).
\end{equation}
We will use this remark with the following result of H\"{o}rmander
(see \cite {HOR1}, tome 4, proposition 25.3.3). Let be $(T,\tau
,x_{0},\xi _{0},y_{0},-\eta _{0})\in \Lambda _{\mathrm{flow}}$,
$\eta _{0}\neq 0$, then near this point there exists, after
perhaps a change of local coordinates in $y$ near $y_{0},$ a
function $S(t,x,\eta )$ such that :
\begin{equation}
\phi (t,x,y,\eta )=S(t,x,\eta )-\left\langle y,\eta \right\rangle
,
\end{equation}
parameterizes $\Lambda _{\mathrm{flow}}$. In particular this
implies that :
\begin{equation*}
\{(t,\partial _{t}S(t,x,\eta ),x,\partial _{x}S(t,x,\eta
),\partial _{\eta }S(t,x,\eta ),-\eta )\}\subset
\Lambda_{\mathrm{flow}} ,
\end{equation*}
and that the function $S$ is a generating function of the flow,
i.e. :
\begin{equation}
\Phi _{t}(\partial _{\eta }S(t,x,\eta ),\eta ) =(x,\partial
_{x}S(t,x,\eta )). \label{Gene}
\end{equation}
Moreover, $S$ satisfies the Hamilton-Jacobi equation :
\begin{equation*}
\left\{
\begin{array}{c}
\partial _{t}S(t,x,\eta )+ p_0(x,\partial
_{x}S(t,x,\eta ))=0, \\
S(0,x,\xi)=\left\langle x,\xi \right\rangle .
\end{array}
\right.
\end{equation*}
Now, we apply this result with $(x_{0},\xi _{0})=(y_{0},\eta
_{0})$, our unique fixed point of the flow on the energy surface
$\Sigma _{E_{c}}$.
\begin{remark}\rm{If $\xi _{0}=0$ we can replace the operator
$P_{h}$ by $e^{\frac{i}{h}\left\langle x,\xi _{1}\right\rangle }P_{h}e^{-\frac{i}{h}%
\left\langle x,\xi _{1}\right\rangle }$ with $\xi _{1}\neq 0.$
This will not change the spectrum since this new operator has the
symbol $p(x,\xi -\xi _{1},h)$ and the critical point is now
$(x_{0},\xi _{1})$ with $\xi _{1}\neq 0$.}
\end{remark}
Consequently, the localized trace $\gamma _{2}(E_{c},h),$ defined
in Eq. (\ref{Trace S2(h)}), can be written for all $N\in
\mathbb{N}$ and modulo $\mathcal{O}(h^{N})$ as :
\begin{equation}
\gamma _{2}(E_{c},h)=\sum\limits_{j<N}(2\pi h)^{-d+j}\int\limits_{\mathbb{%
R\times R}^{2n}}e^{\frac{i}{h}(S(t,x,\xi )-\left\langle x,\xi
\right\rangle +tE_{c})}a_{j}(t,x,\xi )\hat{\varphi}(t)dtdxd\xi .
\label{gamma1 OIF}
\end{equation}
To obtain the right power $-d$ of $h$ occurring here we apply
results of Duistermaat \cite{DUI1} (following here H\"{o}rmander
for the FIO, see \cite {HOR2} tome 4, for example) concerning the
order. An $h$-pseudo-differential operator obtained by Weyl
quantization :
\begin{equation*}
(2\pi h)^{-\frac{N}{2}}\int\limits_{\mathbb{R}^{N}}a(\frac{x+y}{2},\xi )e^{%
\frac{i}{h}\left\langle x-y,\xi \right\rangle }d\xi ,
\end{equation*}
is of order 0 w.r.t. $1/h$. Now since the order of $U_{h}(t)\Theta (P_{h})$ is $-%
\frac{1}{4}$, we find that
\begin{equation}
\psi ^{w}(x,hD_{x})U_{h}(t)\Theta (P_{h})\sim
\sum\limits_{j<N}(2\pi
h)^{-n+j}\int\limits_{\mathbb{R}^{n}}a_{j}(t,x,y,\eta )e^{\frac{i}{h}%
(S(t,x,\eta )-\left\langle y,\eta \right\rangle )}dy.
\label{operateur d'evolution}
\end{equation}
Multiplying Eq. (\ref{operateur d'evolution}) by $\hat{\varphi}%
(t)e^{\frac{i}{h}tE_{c}}$ and passing to the trace we find Eq.
(\ref{gamma1 OIF}) with $d=n$ and where we write again
$a_{j}(t,x,\eta )$ for $a_{j}(t,x,x,\eta )$.

To each element $u_{h}$ of $I(\mathbb{R}^{k},\Lambda )$ we can
associate a principal symbol $e^{\frac{i}{h}S}\sigma
_{\mathrm{princ}}(u_{h})$, where $S$ is a function on $\Lambda $
such that $\xi dx=dS$ on $\Lambda .$ In fact, if
$u_{h}=I(ae^{\frac{i}{h}\varphi })$ then we have $S=S_{\varphi
}=\varphi \circ i_{\varphi }^{-1}$ and $\sigma
_{\mathrm{princ}}(u_{h})$ is a section of $|\Lambda
|^{\frac{1}{2}} \otimes M(\Lambda )$, where $M(\Lambda )$ is the
Maslov vector-bundle of $\Lambda $ and $|\Lambda |^{\frac{1}{2}}$
the bundle of half-densities on $\Lambda$. The half-density of the
propagator $U_{h}(t)$ can be easily expressed in the global
coordinates $(t,y,\eta )$ on $\Lambda _{\mathrm{flow}}.$ If
$p_{1}$ is the sub-principal symbol of $P_{h}$, then this
half-density is given by :
\begin{equation}
\exp (i\int\limits_{0}^{t}p_{1}(\Phi _{s}(y,-\eta ))ds)|dtdyd\eta |^{\frac{1%
}{2}}.\label{demi densite}
\end{equation}
This expression is related to the resolution of the first
transport equation for the propagator, for a proof we refer to
Duistermaat and H\"{o}rmander \cite{D-H}.
\section{Classical dynamic near the equilibrium.}
A critical points of the phase function of (\ref{gamma1 OIF})
satisfies the equations :
\begin{equation*}
\left\{
\begin{array}{c}
E_{c}=-\partial _{t}S(t,x,\xi ), \\
x=\partial _{\xi }S(t,x,\xi ), \\
\xi =\partial _{x}S(t,x,\xi ),
\end{array}
\right. \Leftrightarrow \left\{
\begin{array}{c}
p_0(x,\xi )=E_{c}, \\
\Phi _{t}(x,\xi )=(x,\xi ),
\end{array}
\right.
\end{equation*}
where the right hand side defines a closed trajectory of the flow
inside $\Sigma _{E_{c}}$. Since we are interested in the
contribution of the critical point, we choose a function $\psi \in
C_{0}^{\infty }(T^{\ast }\mathbb{R}^{n})$, with $\psi =1\text{
near }z_{0}$, hence :
\begin{gather*}
\gamma _{2}(E_{c},h) =\frac{1}{2\pi }\mathrm{Tr}\int\limits_{\mathbb{R}}e^{i%
\frac{tE_{c}}{h}}\hat{\varphi}(t)\psi ^{w}(x,hD_{x})\mathrm{exp}(-\frac{i}{h}%
tP_{h})\Theta (P_{h})dt\\
+\frac{1}{2\pi }\mathrm{Tr}\int\limits_{\mathbb{R}}e^{i\frac{tE_{c}}{h}}\hat{%
\varphi}(t)(1-\psi
^{w}(x,hD_{x}))\mathrm{exp}(-\frac{i}{h}tP_{h})\Theta (P_{h})dt.
\end{gather*}
Under the additional hypothesis of having a clean flow, the
asymptotics of the second term is given by the semi-classical
trace formula on a regular level. We also observe that the
contribution of the first term is micro-local. Hence this allows
to introduce local coordinates near $z_{0}$. To separate the
contribution of $z_{0}$ from other closed trajectories we use the
following result on the classical dynamic.
\begin{lemma}\label{dynamic}
For all $T>0$ there exists a neighborhood $U_{T}$ of the critical
point such that $\Phi _{t}(z)\neq z$ for all $z\in U_{T}\backslash
\{z_{0}\}$ and for all $t\in ]\text{-}T,0[\cup ]0,T[$.
\end{lemma}
\noindent \textit{Proof.} Since $z_{0}$ is a degenerate critical
point we have $dH_{p_0}(z_{0})=0$. Hence, for all $\varepsilon >0$
we can find a neighborhood $U$ of $z_{0}$ such that :
\begin{equation*}||dH_{p_0}(z)||\leq \varepsilon, \text{ }\forall z \in
U.
\end{equation*}
By a theorem of Yorke \cite {Yor} we obtain that any closed
trajectory in $U$ has a period $T_{0}$ that satisfies $T_{0}\geq
2\pi \varepsilon^{-1}$. Thus for any $T$ we can choose
$\varepsilon _{T},$ and then $U_{T}$, such that $T_{0}>T$.
\hfill{$\blacksquare$}\medskip

With $\mathrm{supp}(\hat{\varphi})$ compact we can choose $\psi$
such that Lemma \ref{dynamic} holds on $\mathrm{supp}(\psi )$ for
all $t\in \rm{supp}(\hat{\varphi})$. Hence, on the support of
$\psi $ there is two contributions :\\
1) Points $(t,x,\xi )=(0,x,\xi )$ for $(x,\xi )\in \Sigma
_{E_{c}}$.\\
2) Points $(t,x,\xi )=(t,z_{0})$ for $t\in
\rm{supp}(\hat{\varphi})$.\newline The first contribution is
non-singular for $(x,\xi )\neq z_{0}$ and can be treated, again,
by the regular trace formula. Now, we restrict our attention to
the second contribution and, since $z_0$ is totally degenerate, we
obtain :
\begin{equation}
d\Phi _{t}(z_{0})=\mathrm{exp}(0)=\textrm{Id},\text{ }\forall t.
\end{equation}
The next homogeneous components of the flow are given, with
$(H_{2})$, by :
\begin{equation}d^{j}\Phi
_{t}(z_{0})=0,\text{ }\forall t,\text{ }\forall j\in \{2,..,k-2\}.
\end{equation}
To obtain the next non-zero term of the Taylor expansion of the
flow, we will use the following technical result :
\begin{lemma}
\label{TheoFormule de récurence du flot}Let be $z_{0}$ an
equilibrium of the $C^{\infty}$ vector field $X$ and $\Phi _{t}$
the flow of $X$. Then for all $m\in \mathbb{N}^{\ast }$, there
exists a polynomial map $P_{m}$, vector valued and of degree at
most $m$, such that :
\begin{equation*}
d^{m}\Phi _{t}(z_{0})(z^{m})=d\Phi
_{t}(z_{0})\int\limits_{0}^{t}d\Phi _{-s}(z_{0})P_{m}(d\Phi
_{s}(z_{0})(z),...,d^{m-1}\Phi _{s}(z_{0})(z^{m-1}))ds.
\end{equation*}
\end{lemma}
For a proof we refer to \cite{Cam0} or \cite{Cam}.\medskip

Since $d\Phi_{t}(z_{0})=\mathrm{Id}$, for all $t$, the first
non-zero term of the Taylor expansion of the flow is given by :
\begin{equation}
d^{k-1}\Phi_{t}(z_{0})(z^{k-1})=\int\limits_{0}^{t}
d^{k-1}H_{p_0}(z_{0})(z^{k-1})ds=td^{k-1}H_{\mathfrak{p}_{k}}(z_{0})(z^{k-1}).
\label{derive ordre k-1 du flot}
\end{equation}
Where the last identity is obtained using that
$d^{k-1}H_{p_0}(z_0)=d^{k-1}H_{\mathfrak{p}_k}(z_0)$. Moreover,
with $d^{2}p(z_0)=0$, for the next term Lemma \ref{TheoFormule de
récurence du flot} gives :
\begin{equation*}
d^{k}\Phi_{t}(z_{0})(z^{k})=\int\limits_{0}^{t}d^{k}H_{p_0}(z_{0})(z^{k})ds
=td^{k}H_{\mathfrak{p}_{k+1}}(z_{0})(z^{k}).
\end{equation*}
\begin{remark}\label{struc flow}
\rm{For higher order derivatives $d^{j}\Phi_{t}(z_0)$, with $j>k$,
there is two different kind of terms, namely :
\begin{equation*}
\int\limits_{0}^{t}d^{j}H_{p_0}(z_{0})(z^{j})ds
=td^{j}H_{p_0}(z_{0})(z^{j})=\mathcal{O}(||z||^{j}),
\end{equation*}
and terms involving powers of $t$, for example we have :
\begin{gather*}
\int\limits_{0}^{t}d^{j+2-k}H_{p_0}(z_{0})(z^{j+1-k},d^{k-1}\Phi_{t}(z_{0})(z^{k-1}))\\
=\frac{t^{2}}{2}d^{j+2-k}H_{p_0}(z_{0})(z^{j+1-k},d^{k-1}H_{\mathfrak{p}_{k}}(z_{0})(z^{k-1})).
\end{gather*}
This term is simultaneously $\mathcal{O}(t^2)$ and
$\mathcal{O}(||z||^j)$ near $(0,z_0)$. A similar result holds for
other terms, which are $\mathcal{O}(t^d)$ for $d\geq2$, by an easy
recurrence.}
\end{remark}
\begin{lemma}\label{structure phase} Near $z_0$, here supposed to be 0 to simplify, we
have :
\begin{equation}
S(t,x,\xi )-\left\langle x,\xi \right\rangle
+tE_{c}=-t(\mathfrak{p}_{k}(x,\xi )+R_{k+1}(x,\xi
)+tG_{k+1}(t,x,\xi )), \label{forme phase}
\end{equation}
where $R_{k+1}(x,\xi ) =\mathcal{O}(||(x,\xi )||^{k+1})$ and
$G_{k+1}(t,x,\xi)=\mathcal{O}(||(x,\xi )||^{k+1})$, uniformly with
respect to $t$.
\end{lemma}
\textit{Proof.} By Taylor we obtain :
\begin{equation*}
\Phi _{t}(x,\xi ) =(x,\xi )+\frac{1}{(k-1)!}d^{k-1}\Phi
_{t}(0)(z^{k-1})+\mathcal{O}(||z||^{k}).
\end{equation*}
We search our local generating function as :
\begin{equation*}
S(t,x,\xi)=-tE_c+\left\langle x,\xi
\right\rangle+\sum\limits_{j=3}^N
S(t,x,\xi)+\mathcal{O}(||(x,\xi)||^{N+1}),
\end{equation*}
where the $S_j$ are time dependant and homogeneous of degree $j$
w.r.t. $(x,\xi)$. With the implicit relation
$\Phi_t(\partial_{\xi}S(t,x,\xi),\xi)=(x,\partial_{x}S(t,x,\xi))$,
we have :
\begin{equation*}
S(t,x,\xi)=-tE_c+\left\langle x,\xi \right\rangle+S_{k}(t,x,\xi)+
\mathcal{O}(||(x,\xi )||^{k+1}),
\end{equation*}
where $S_k$ is homogeneous of degree $k$ w.r.t. $(x,\xi)$. If $J$
is the matrix of the usual symplectic form, comparing terms of
same degree gives :
\begin{equation*}
J\nabla S_k(t,x,\xi)=\frac{1}{(k-1)!}d^{k-1}\Phi
_{t}(0)((x,\xi)^{k-1}),
\end{equation*}
By homogeneity and with Eq. (\ref{derive ordre k-1 du flot}) we
obtain :
\begin{gather*}
S_k(t,x,\xi)=\frac{1}{k!}\left\langle
(x,\xi),tJd^{k-1}H_{\mathfrak{p}_k}(x,\xi)^{k-1}
\right\rangle=-t\mathfrak{p}_{k}(x,\xi),\\
S(t,x,\xi)=-tE_c+\left\langle x,\xi
\right\rangle-t\mathfrak{p}_{k}(x,\xi) +\mathcal{O}(||(x,\xi
)||^{k+1}).
\end{gather*}
As concern the remainder, we first observe that $S(0,x,\xi)=
\left\langle x,\xi \right\rangle $. Hence, we can write :
\begin{equation*}
S(t,x,\xi)-\left\langle x,\xi \right\rangle=tF(t,x,\xi),
\end{equation*}
with $F$ smooth in a neighborhood of $(x,\xi)=0$. Now, the
Hamilton-Jacobi equation imposes that $F(0,x,\xi)=-p_0(x,\xi)$ and
we obtain :
\begin{equation*}
R_{k+1}(x,\xi)=p_0(x,\xi)-E_c-\mathfrak{p}_k(x,\xi)=\mathcal{O}(||(x,\xi)||^{k+1}).
\end{equation*}
Finally, the time dependant remainder can be written :
\begin{equation*}
S(t,x,\xi)-S(0,x,\xi)-t\partial_t S(0,x,\xi)=\mathcal{O}(t^2),
\end{equation*}
since by construction this term is of order
$\mathcal{O}(||(x,\xi)||^{k+1})$ we get the desired result when
$t$ is in a compact subset of $\mathbb{R}$. $\hfill{\blacksquare}$
\section{Normal forms of the phase function.}
Since the contribution we study is local, we can work with some
coordinates and identify locally $T^{\ast }\mathbb{R}^{n}$ with
$\mathbb{R}^{2n}$ near the critical point. We define :
\begin{equation}
\Psi(t,z)=\Psi(t,x,\xi)=S(t,x,\xi )-\left\langle x,\xi
\right\rangle+tE_{c},\text{ }z=(x,\xi)\in \mathbb{R}^{2n}.
\label{defphase}
\end{equation}
\begin{lemma}\label{FN1}
If $P_h$ satisfies conditions $(H_{2})$ and $(H_{4})$ then, in a
neighborhood of $(t,z)=(0,z_0)$, there exists local coordinates
$\chi$ such that :
\begin{gather*}
\Psi(t,z) \simeq -\chi _{0}\chi_{1}^{k},\text{ in all
directions where }p_{k}(\theta )>0,\\
\Psi(t,z) \simeq + \chi _{0}\chi_{1}^{k},\text{ in all
directions where }p_{k}(\theta )<0,\\
 \Psi(t,z) \simeq \chi
_{0}\chi _{1}^{k}\chi_{2},\text{near any point where
}p_{k}(\theta)= 0.
\end{gather*}
\end{lemma}
\noindent\textit{Proof.} We can here assume that $z_{0}$ is the
origin and we use polar coordinates $z=(r,\theta ),$ $\theta \in
\mathbb{S}^{2n-1}(\mathbb{R})$. With Lemma \ref{structure phase},
near the critical point we have :
\begin{equation*}
\Psi(t,z) \simeq -tr^{k}(\mathfrak{p}_{k}(\theta )+rR_{k+1}(\theta
)+tG_{k+1}(t,r\theta )),
\end{equation*}
where $\mathfrak{p}_{k}(\theta )$ is the restriction of
$\mathfrak{p}_{k}$ on
$\mathbb{S}^{2n-1}$.\\
If $\mathfrak{p}_{k}(\theta_{0})\neq 0$ we define new coordinates
:
\begin{gather*}
(\chi _{0},\chi _{2},...,\chi _{2n})(t,r,\theta) =(t,\theta _{1},...,\theta _{2n-1}), \\
\chi _{1}(t,r,\theta) = r|\mathfrak{p}_{k}(\theta
)+rR_{k+1}(\theta )+tG_{k+1}(t,r\theta )|^{\frac{1}{k}}.
\end{gather*}
In these coordinates the phase becomes $-\chi _{0}\chi _{1}^{k}$
if $\mathfrak{p}_{k}(\theta_0)$ is positive (resp. $\chi _{0}\chi
_{1}^{k}$ for a negative value). Near $\theta_0$, we have :
\begin{equation*}
\frac{\partial \chi _{1}}{\partial r} (t,0,\theta) =
|\mathfrak{p}_{k}(\theta )|^{\frac{1}{k}}\neq 0, \text{ } \forall
t,
\end{equation*}
hence, the corresponding Jacobian satisfies $|J\chi |(t,0,\theta
)=|\mathfrak{p}_{k}(\theta )|^{\frac{1}{k}}\neq 0$.

Now, let $\theta_0$ be such that $\mathfrak{p}_{k}(\theta_{0})=
0$. Up to a permutation, we can suppose that
$\partial_{\theta_{1}}\mathfrak{p}_{k}(\theta_{0})\neq 0$. We
choose here the new coordinates :
\begin{gather*}
(\chi _{0},\chi _{1},\chi
_{3},...,\chi_{2n})(t,r,\theta)=(t,r,\theta
_{2},...,\theta _{2n-1}), \\
\chi_{2}(t,r,\theta)=\mathfrak{p}_{k}(\theta)+rR_{k+1}(\theta)+
tG_{k+1}(t,r\theta).
\end{gather*}
Since we have $|J\chi
|(t,0,\theta_0)=|\partial_{\theta_{1}}\mathfrak{p}_{k}(\theta_{0})|\neq
0$, lemma follows. \hfill{$\blacksquare$}

In order to use these normal forms we introduce an adapted
partition of unity on $\mathbb{S}^{2n-1}$. We choose functions
$\Omega_{j}(\theta)$, with compact supports, such that :
\begin{equation*}
\{\theta\in \mathbb{S}^{2n-1}\text{ / }
\mathfrak{p}_{k}(\theta)=0\} \subset \bigcup\limits_{j}
\textrm{supp}(\Omega_{j}),
\end{equation*}
so that normal forms of Lemma \ref{FN1} exist inside
$\textrm{supp}(\Omega_{i})$. Since $\mathbb{S}^{2n-1}$ is compact
this set of functions can be chosen finite and we obtain a
partition of unity with $\Omega_{0}=1-\sum \Omega_{i}$. The
support of $\Omega_{0}$ might be not connected and we define
$\Omega_{0}^{+}$, with $\mathfrak{p}_k(\theta)>0$ on
$\mathrm{supp}(\Omega_{0}^{+})$, and similarly we define
$\Omega_{0}^{-}$ where $\mathfrak{p}_k<0$,
so that $\Omega_{0}=\Omega_{0}^{+}+\Omega_{0}^{-}$.\\
If we accordingly split up our oscillatory-integral we obtain :
\begin{gather*}
I_{+}(\lambda)=\int\limits_{\mathbb{R\times R}_{+}\times
\mathbb{S}^{2n-1}}e^{\frac{i}{h}\Psi (t,r,\theta
)}\Omega_{0}^{+}(\theta) a(t,r\theta )r^{2n-1}dtdrd\theta\\
=\int\limits_{\mathbb{R\times R}_{+}} e^{-\frac{i}{h}\chi _{0}\chi
_{1}^{k}}A_{0}^{+}(\chi _{0},\chi _{1})d\chi _{0}d\chi _{1},
\end{gather*}
for the directions where $\mathfrak{p}_k(\theta)>0$ and :
\begin{gather*}
I_{-}(\lambda)=\int\limits_{\mathbb{R\times R}_{+}\times
\mathbb{S}^{2n-1}}e^{\frac{i}{h}\Psi (t,r,\theta
)}\Omega_{0}^{-}(\theta) a(t,r\theta )r^{2n-1}dtdrd\theta
\\=\int\limits_{\mathbb{R\times R}_{+}} e^{\frac{i}{h}\chi _{0}\chi
_{1}^{k}}A_{0}^{-}(\chi _{0},\chi _{1})d\chi _{0}d\chi _{1},
\end{gather*}
for the directions where $\mathfrak{p}_k(\theta)<0$. Similarly,
the contribution of the neighborhood of the set
$\{\mathfrak{p}_k(\theta)=0\}$ is given by :
\begin{gather*}
I_{j}(\lambda)=\int\limits_{\mathbb{R\times R}_{+}\times
\mathbb{S}^{2n-1}}e^{\frac{i}{h}\Psi (t,r,\theta
)}\Omega_{j}(\theta) a(t,r\theta )r^{2n-1}dtdrd\theta \\
=\int\limits_{\mathbb{R}\times \mathbb{R}_{+}\times\mathbb{R}}
e^{\frac{i}{h}\chi _{0}\chi _{1}^{k}\chi_{2}}A_{j}(\chi _{0},\chi
_{1},\chi_{2})d\chi _{0}d\chi _{1}d\chi_{2}.
\end{gather*}
The associated new amplitudes are respectively given by :
\begin{gather}
A_{0}^{\pm}(\chi _{0},\chi _{1})=\int \chi ^{\ast
}(\Omega_{0}^{\pm}(\theta)a(t,r\theta
)r^{2n-1}|J\chi |)d\chi _{2}...d\chi _{2n}, \label{ampli1}\\
A_{j}(\chi _{0},\chi _{1},\chi_{2})=\int \chi ^{\ast
}(\Omega_{j}(\theta)a(t,r\theta )r^{2n-1}|J\chi |)d\chi
_{3}...d\chi _{2n}. \label{ampli2}
\end{gather}
\begin{remark}\label{degres amplitude}
\rm{Since $\chi _{1}(t,r,\theta ) =r|\mathfrak{p}_{k}(\theta
)+rR_{k+1}(\theta )+tG_{k+1}(t,r\theta )|^{\frac{1}{k}}$, our new
amplitude satisfies $A_{0}^{\pm}(\chi _{0},\chi
_{1})=\mathcal{O}(\chi _{1}^{2n-1})$, near $\chi _{1}=0$. A
similar argument shows that $A_{i}(\chi _{0},\chi
_{1},\chi_{2})=\mathcal{O}(\chi _{1}^{2n-1})$, near $\chi _{1}=0$.
These facts will play a major role in Lemmas \ref{Theo IO 1ere
carte}, \ref{Theo IO 2eme carte} and \ref{Theo IO 2eme with log}
below. }
\end{remark}
We end this section with lemmas on asymptotics of oscillatory
integrals with phases as in Lemma \ref{FN1}.
\begin{lemma}
\label{Theo IO 1ere carte}There exists a sequence $(c_{j})_{j}$ of
distributions, whose support is contained in the set $\{\chi
_{1}=0\},$ such that for all function $a\in C_{0}^{\infty }(\mathbb{R}%
_{+}$ $\times \mathbb{R})$ :
\begin{equation}
\int\limits_{0}^{\infty }(\int\limits_{\mathbb{R}}e^{i\lambda \chi
_{0}\chi _{1}^{k}}a(\chi _{0},\chi _{1})d\chi _{0})d\chi _{1}\sim
\sum\limits_{j=0}^{\infty }\lambda ^{-\frac{j+1}{k}}c_{j}(a),
\end{equation}
asymptotically for $\lambda \rightarrow \infty ,$ where :
\begin{equation*}
c_{j}=\frac{1}{k}\frac{1}{j!}(\mathcal{F}(x_{-}^{\frac{j+1-k}{k}})(\chi_0)
\otimes\delta _{0}^{(j)}(\chi_1)),\text{ }x_{-}=\max (-x,0).
\end{equation*}
\end{lemma}
We refer to \cite{Cam1} for a proof of this lemma.
$\hfill{\blacksquare}$
\begin{remark}
\rm{A similar result holds for a phase $-\chi _{0}\chi _{1}^{k}$
if we replace terms $x_{-}$ by $x_{+}$ in Lemma \ref{Theo IO 1ere
carte}.}
\end{remark}
\begin{lemma}\label{Theo IO 2eme carte}
If $k>2n$ we have :
\begin{equation*}
\int\limits_{0}^{\infty }(\int\limits_{\mathbb{R}^2}e^{i\lambda
\chi _{0}\chi _{1}^{k}\chi_2}a(\chi _{0},\chi _{1},\chi_2)d\chi
_{0}d\chi_2)\chi _{1}^{2n-1}d\chi _{1}= \lambda
^{-\frac{2n}{k}}d(a)+\mathcal{O}(\lambda
^{-\frac{2n+1}{k}}\mathrm{log}^2(\lambda)),
\end{equation*}
where the leading coefficient is given by :
\begin{equation} \label{toporder1}
d(a)=\frac{1}{k}\Gamma(\frac{2n}{k}) \int\limits_{\mathbb{R}^2}
|\chi_0 \chi_2 |^{-\frac{2n}{k}}\mathrm{exp}(i\frac{\pi
n}{k}\mathrm{sign}(\chi_0\chi_2)) a(\chi_0,0,\chi_2)d\chi_0
d\chi_2.
\end{equation}
\end{lemma}
\noindent\textit{Proof.} We use the Bernstein-Sato polynomial, see
\cite{WON} for a detailed construction. We use variables $(t,r,v)$
instead of $(\chi_0,\chi_1,\chi_2)$ and since $\chi_1 \geq 0$ for
$tv\geq 0$ we can write :
\begin{gather}
\frac{\partial^{2}}{\partial t\partial v}
\frac{\partial^{k}}{\partial r^{k}}
((tvr^{k})^{1-z}r^{2n-1})=b_k(z)(tvr^{k})^{-z}r^{2n-1}, \label{trick}\\
b_k(z)=(1-z)^{2}\prod\limits_{j=1}^{k}(j-kz+2n-1)\label{trick2}.
\end{gather}
With the classical representation :
\begin{gather*} \int\limits_{0}^{\infty
}\int\limits_{\{tv \geq 0\}}e^{i\lambda
tvr^{k}}a(t,r,v)dtdvdr\\
=\frac{1}{2i\pi }\int\limits_{\gamma}
e^{i\frac{\pi z}{2}}\Gamma (z) \lambda ^{-z}
(\int\limits_{0}^{\infty }\int\limits_{\{tv \geq 0\}}
(tvr^{k})^{-z}a(t,r,v)dtdvdr)dz,
\end{gather*}
where $\gamma=]c-i\infty ,c+i\infty [$ and
$\textrm{Re}(c)<k^{-1}$, we can compute the asymptotic expansion
with the residue method by pushing of the complex path of
integration $\gamma$ to the right. Note that, when the phase is
negative, we have :
\begin{equation}
\frac{1}{2i\pi }\int\limits_{\gamma} e^{-i\frac{\pi z}{2}}\Gamma
(z) \lambda ^{-z} ( \int\limits_{0}^{\infty }\int\limits_{\{tv
\leq 0\}} (|tv|r^{k})^{-z}a(t,r,v)dtdvdr)dz.\label{phase negative}
\end{equation}
With Eq. (\ref{trick}) and (\ref{trick2}) the meromorphic
extension is given by :
\begin{gather*}
\int\limits_{0}^{\infty }\int\limits_{\{tv \geq 0\}}
(tr^{k}v)^{-z}r^{2n-1}a(t,r,v)dtdrdv\\
=\frac{(-1)^{k}}{b_k(z)}\int\limits_{0}^{\infty }\int\limits_{\{tv
\geq 0\}} (tvr^{k})^{1-z}r^{2n-1}\frac{\partial^{2}}{\partial
t\partial v} \frac{\partial^{k}}{\partial r^{k}}a(t,r,v)dtdvdr.
\end{gather*}
Under our assumptions, the poles are $z=1$ as a double pole and :
\begin{equation*}
z=\frac{j+2n-1}{k},\text{ } j\in[1,...,k].
\end{equation*}
Since $k>2n$, the first pole is $\frac{2n}{k}\notin \mathbb{Z}$.
The residue in this pole is :
\begin{gather*}
c_{k,n}\lambda^{-\frac{2n}{k}}\mathrm{exp}(i\pi\frac{n}{k})\Gamma
(\frac{2n}{k})\int\limits_{\{tv>0\}} (tv)^{1-\frac{2n}{k}}
r^{k-1}\frac{\partial^{2}}{\partial t\partial v}
\frac{\partial^{k}}{\partial r^{k}}a(t,r,v)dtdvdr,\\
c_{k,n}=(-1)^k \lim_{z\rightarrow \frac{2n}{k}} \frac
{(z-\frac{2n}{k})}{b_k(z)}=(-1)^{k+1} \frac{k}{(k-2n)^2\Gamma
(k)}.
\end{gather*}
Now, with the following relations :
\begin{gather*}
\int\limits_{t,v>0}
(tv)^{1-\frac{2n}{k}}\frac{\partial^{2}a}{\partial t\partial
v}(t,v)dtdv=(1-\frac{2n}{k})^2\int\limits_{t,v>0}
(tv)^{-\frac{2n}{k}}a(t,v)dtdv, \\
(-1)^{k}\int\limits_{r>0} r^{k-1}\frac{\partial^{k}a}{\partial
r^k}dr=(k-1)!a(0),
\end{gather*}
and also using that :
\begin{equation*}
\frac{k}{(k-2n)^2\Gamma (k)} (1-\frac{2n}{k})^2 (k-1)! \Gamma
(\frac{2n}{k})=\frac{\Gamma (\frac{2n}{k})}{k},
\end{equation*}
we find that the first residue is given by :
\begin{equation*}
\frac{\mathrm{exp}(i\pi\frac{n}{k})}{k}\lambda^{-\frac{2n}{k}}\Gamma
(\frac{2n}{k}) \int\limits_{\{t,v\geq 0\}} (tv)^{-\frac{2n}{k}}
a(t,0,v)dtdv.
\end{equation*}
The summation  over the other quadrants gives the result with Eq.
(\ref{phase negative}). Finally, the remainder is of order
$\mathcal{O}(\lambda ^{-\frac{2n+1}{k}}\mathrm{log}^2(\lambda))$
if $2n+1=k$ and of order  $\mathcal{O}(\lambda
^{-\frac{2n+1}{k}})$ otherwise.\hfill{$\blacksquare$}
\begin{remark}
\rm{The preceding result holds again, in a weaker sense, if $k$ is
not a divisor of $2n$, as can be shown by iterating the process of
meromorphic extension above and using the fact that our amplitude
is $\mathcal{O}(\chi_{1}^{2n-1})$ near $\chi_1 =0$.\newline%
Now, if $k$ divide $2n$ then terms
$\lambda^{-j}\mathrm{log}^d(\lambda)$, with $d=1,2$, can appear
for the leading term, since we obtain poles of order 3 for the
first non-zero residue.}
\end{remark}
\textbf{Example.}\\
If we consider Gaussian amplitude, we obtain easily :
\begin{equation*}
I(\lambda,k,n)=\int\limits_{\mathbb{R}_{+} \times \mathbb{R}^2}
e^{i\lambda tr^k v} e^{-(t^2+v^2+r^2)}r^{2n-1} drdtdv  =2\pi
\int\limits_{0}^{\infty} \frac{r^{2n-1}e^{-r^2}}{\sqrt{4+r^{2k}
\lambda^2}}dr.
\end{equation*}
The choice of $k=5$ and $n=2$ leads to :
\begin{equation*}
I(\lambda,5,2)\sim 2\pi \frac{\Gamma(\frac{2}{5}) \Gamma
(\frac{11}{10})}{2^{\frac{1}{5}} \sqrt{\pi}}
\lambda^{-\frac{4}{5}}+\mathcal{O}(\lambda^{-1}).
\end{equation*}
This extra Gamma-factor comes from the identity :
\begin{equation*}
\int\limits_{\mathbb{R}^2} |tv|^{-\frac{4}{5}} e^{-t^2} e^{-v^2}
dtdv =\Gamma ( \frac{1}{10})^2.
\end{equation*}
\begin{lemma}\label{Theo IO 2eme with log}
When $k=2n$ we have the particular result :
\begin{equation*}
\int\limits_{0}^{\infty }(\int\limits_{\mathbb{R}^2}e^{i\lambda
\chi _{0}\chi _{1}^{2n}\chi_2}a(\chi _{0},\chi _{1},\chi_2)d\chi
_{0}d\chi_2)\chi _{1}^{2n-1}d\chi _{1}= \lambda
^{-1}\mathrm{log}(\lambda)d(a)+\mathcal{O}(\lambda ^{-1}),
\end{equation*}
with the leading coefficient :
\begin{equation}
d(a)=\frac{\pi}{n}a(0,0,0).
\end{equation}
\end{lemma}
\noindent\textit{Proof.} We use again the Bernstein-Sato
polynomial method. We use variables $(t,r,v)$ instead of
$(\chi_0,\chi_1,\chi_2)$ and we define :
\begin{equation*}
J_{+}(\lambda)=\int\limits_{0}^{\infty}\int\limits_{\{tv>0\}
}e^{i\lambda tvr^{2n}}a(t,r,v)r^{2n-1}dtdrdv,
\end{equation*}
similarly, we define $J_{-}(\lambda )$ on the set $\{tv <0\}$.
Then we can write :
\begin{equation*}
J_{+}(\lambda)=\frac{1}{2i\pi}\int\limits_{\gamma} e^{i\frac{\pi
z}{2}} \Gamma (z) \lambda^{-z}
\int\limits_{0}^{\infty}\int\limits_{\{tv>0\}} (tv)^{-z} r^{-2nz}
a(t,r,v)r^{2n-1}dtdrdv,
\end{equation*}
where $\gamma=]c-i\infty,c+i\infty[$, $c<(2n)^{-1}$. Similar
computations as in the proof of Lemma (\ref{Theo IO 2eme carte})
show that the associated Bernstein-Sato polynomial is :
\begin{equation*}
b_{2n}(z)=(1-z)^{2}\prod\limits_{j=1}^{2n}(j-2nz+2n-1),
\end{equation*}
under our assumptions the first pole is $z=1$ and is of order 3.
We define two holomorphic functions, near $z=1$, via :
\begin{equation*}
G^{\pm}_n(z)=\frac{(z-1)^3}{b_{2n}(z)} e^{ \pm i\frac{\pi
z}{2}}\Gamma (z).
\end{equation*}
Hence, when using the residue method, the first terms of the
asymptotic expansion of $J_{+}(\lambda)$ are given by the formula
:
\begin{equation*}
\frac{1}{2}\lim\limits_{z\rightarrow 1}\left (
\frac{\partial^2}{\partial z^2} (G^{+}_{n}(z) \lambda^{-z}
\int\limits_{0}^{\infty}\int\limits_{\{tv>0\}} (tv)^{1-z}
r^{2n-2nz}r^{2n-1} \frac{\partial^2}{\partial t
\partial v} \frac{\partial^{2n}}{\partial r ^{2n}}a(t,r,v)dtdrdv )\right ),
\end{equation*}
and we have a similar formula for $J_{-}(\lambda)$ if we use
$G^{-}_{n}(z)$ and integration over $\{tv<0\}$. For all
holomorphic application $f$ and all $\lambda
>0$, we have :
\begin{equation}
\frac{\partial^2}{\partial z^2}
(f(z)\lambda^{-z})=\frac{\partial^2f}{\partial
z^2}(z)\lambda^{-z}-2\log(\lambda)\frac{\partial f}{\partial
z}(z)\lambda^{-z}+\log^2(\lambda)f(z)\lambda^{-z}.
\label{splitting}
\end{equation}
The term involving $\log^2(\lambda)$ is computed as in Lemma
\ref{Theo IO 2eme carte} and we have :
\begin{equation*}
\lim\limits_{z\rightarrow 1} G^{+}_n(z)=- \frac{i}{(2n)!}.
\end{equation*}
The distributional factor is here given by :
\begin{equation*}
\int\limits_{0}^{\infty}\int\limits_{\{tv>0\}} r^{2n-1}
\frac{\partial^2}{\partial t
\partial v} \frac{\partial^{2n}}{\partial r
^{2n}}a(t,r,v)dtdrdv=2(2n-1)! a(0,0,0).
\end{equation*}
Hence the contribution is given by :
\begin{equation*}
-\log^2(\lambda)\lambda^{-1} \frac{i}{n}a(0,0,0),
\end{equation*}
but when the phase is negative we have :
\begin{equation*}
\lim\limits_{z\rightarrow 1} G^{-}_n(z)= + \frac{i}{(2n)!},
\end{equation*}
consequently, by summation, there is no term associated to
$\log(\lambda)^2$.\medskip\\
Now, we compute the main term associated to the $\log
(\lambda)$.\\
For the same reason as previously all terms obtained by derivation
of $\Gamma(z)$ and $(z-1)^3/b_{2n}(z)$ will give a zero
contribution by summation. This will not be true for terms
obtained by derivation of the exponential. As concerns derivation
of the meromorphic distributions, we obtain respectively :
\begin{gather*}
- \frac{i}{(2n)!} \int\limits_{0}^{\infty}\int\limits_{\{tv>0\}}
\log (tv r^{2n}) r^{2n-1} \frac{\partial^2}{\partial t
\partial v} \frac{\partial^{2n}}{\partial r
^{2n}}a(t,r,v)dtdrdv,\\
+ \frac{i}{(2n)!} \int\limits_{0}^{\infty}\int\limits_{\{tv<0\}}
\log (|tv| r^{2n}) r^{2n-1} \frac{\partial^2}{\partial t
\partial v} \frac{\partial^{2n}}{\partial r
^{2n}}a(t,r,v)dtdrdv.
\end{gather*}
By summation and using that $\log
(|tv|r^{2n})=\log(|tv|)+\log(r^{2n})$ we easily obtain that the
associated factor is given by :
\begin{equation*}
-\frac{i}{2n} \left (\int\limits_{\{tv>0\}} \log (tv)
\frac{\partial^2}{\partial t
\partial v} a(t,0,v)dtdv-  \int\limits_{\{tv<0\}} \log (|tv|)
\frac{\partial^2}{\partial t
\partial v} a(t,0,v)dtdv \right ).
\end{equation*}
A new splitting of the logarithms and integrations by parts show
that the associated contribution vanish.

It remains now to compute the term associated to the derivation of
the exponential. An easy computation shows that the associated
contribution is :
\begin{equation*}
\log(\lambda)\lambda^{-1}
(i\frac{\pi}{2})\frac{e^{i\frac{\pi}{2}}}{(2n)!}
 \int\limits_{0}^{\infty}\int\limits_{\{tv>0\}}
r^{2n-1} \frac{\partial^2}{\partial t
\partial v} \frac{\partial^{2n}}{\partial r ^{2n}}a(t,r,v)dtdrdv,
\end{equation*}
by integrations by parts, we finally obtain that the contribution
is given by :
\begin{equation*}
\log(\lambda)\lambda^{-1}
(i\frac{\pi}{2})\frac{e^{i\frac{\pi}{2}}}{2n}
\int\limits_{\{tv>0\}} \frac{\partial^2}{\partial t
\partial v} a(t,0,v)dtdv=-\frac{\pi}{2n} \log(\lambda)\lambda^{-1}
a(0,0,0).
\end{equation*}
With Eq. (\ref{splitting}) we have :
\begin{equation*}
J_{+}(\lambda)=\frac{\pi}{2n} \log(\lambda)\lambda^{-1} a(0,0,0)
+\mathcal{O}(\lambda^{-1}).
\end{equation*}
A totally similar calculation gives :
\begin{equation*}
J_{-}(\lambda)=\frac{\pi}{2n} \log(\lambda)\lambda^{-1} a(0,0,0)
+\mathcal{O}(\lambda^{-1}),
\end{equation*}
and the result holds by summation of $J_{-}(\lambda)$ and
$J_{+}(\lambda)$.$\hfill{\blacksquare}$\medskip\\
\textbf{Example.}\\
We use again an amplitude that is a product of Gaussian-functions,
we have :
\begin{equation*}
I(\lambda,k,n)=\int\limits_{\mathbb{R}_{+} \times \mathbb{R}^2}
e^{i\lambda tr^k v} e^{-(t^2+v^2+r^2)}r^{2n-1} drdtdv  =2\pi
\int\limits_{0}^{\infty} \frac{r^{2n-1}e^{-r^2}}{\sqrt{4+r^{2k}
\lambda^2}}dr.
\end{equation*}
Hence, for $k=2$ and $n=1$ we obtain :
\begin{equation*}
I(\lambda,2,1)=\frac{1}{4\lambda} (-2 \pi^2
Y_{0}(\frac{2}{\lambda}) + (2\pi-1 ) J_{0}( \frac{2}{\lambda})
\log(\frac{1}{\lambda^2}) + \pi\mathrm{H}_0 (\frac{2}{\lambda})),
\end{equation*}
where  $J_v (x)$, $Y_v(x)$  are respectively the standard Bessel
functions of first and second kind. Also, $ \mathrm{H}_v (x)$ is
the standard Struve function defined by :
\begin{equation*}
z^2 y''(z)+zy'(z)+(z^2 -v^2)y(z) = \frac{2}{\pi} \frac{z^{v+1}}{(2
v - 1)!!}.
\end{equation*}
From classical properties of these special functions, we obtain :
\begin{equation}
I(\lambda,2,1) \sim (\frac{ \pi \log (\lambda)-\gamma\pi
}{\lambda}) +\mathcal{O}(\lambda^{-2}).
\end{equation}
Here, $\gamma$ is Euler's constant and is obtained by derivation
of the $\Gamma (z)$ factor in the formula that gives meromorphic
extensions of our distributions. $\hfill{\blacksquare}$
\section{Proof of the main results.}
\textbf{Directions where $\mathfrak{p}_k(\theta)\neq 0$.}\\
Following step by step the proof of Lemma 4 of \cite{Cam1} we
obtain that the first non-zero coefficient is obtained for
$l=2n-1$ (see Remark \ref{degres amplitude}) and is given by
\begin{equation*}
\frac{1}{k}\frac{1}{(2n-1)!}\left\langle
\mathcal{F}(x_{+}^{\frac{2n-k}{k}})\otimes \delta
_{0}^{(2n-1)},A_{0}^{+}(\chi _{0},\chi _{1})\right\rangle
=\frac{1}{k}\int \mathcal{F}(x_{+}^{\frac{2n-k}{k}})
(\chi_{0})\tilde{A}_{0}^{+}(\chi _{0},0)d\chi _{0}.
\end{equation*}
Since by construction :
\begin{equation}
\tilde{A}_{0}^{+}(\chi
_{0},0)=\int\limits_{\mathbb{S}^{2n-1}}a(\chi
_{0},0)\Omega_{0}^{+}(\theta)|\mathfrak{p}_{k}(\theta
)|^{-\frac{2n}{k}}d\theta ,
\end{equation}
we obtain that the local contribution, associated to
$\mathrm{supp}(\Omega_{0}^{+})$, is :
\begin{equation}
\frac{1}{k}\left\langle (\mathcal{F}(x_{+}^{\frac{2n-k}{k}})
(\chi_{0}),a(\chi _{0},0)\right\rangle
\int\limits_{\mathbb{S}^{2n-1}}\Omega_{0}^{+}(\theta)|\mathfrak{p}_{k}(\theta
)|^{-\frac{2n}{k}}d\theta.
\end{equation}
A similar computation gives the contribution of
$\mathrm{supp}(\Omega_{0}^{-})$, via :
\begin{equation}
\frac{1}{k}\left\langle (\mathcal{F}(x_{-}^{\frac{2n-k}{k}})
(\chi_{0}),a(\chi _{0},0)\right\rangle
\int\limits_{\mathbb{S}^{2n-1}}\Omega_{0}^{-}(\theta)|\mathfrak{p}_{k}(\theta
)|^{-\frac{2n}{k}}d\theta.
\end{equation}
Now, since $a(t,0)=\hat{\varphi}(t)\exp(itp_1(z_0))$, cf. Eq.
(\ref{demi densite}), the contributions of the directions where
$\mathfrak{p}_k(\theta)\neq 0$ are given, respectively, by :
\begin{equation} \label{positive}
I_{+}(\lambda) \sim \frac{1}{k} \lambda^{-\frac{2n}{k}}
\left\langle |t|_{+}^{\frac{2n-k}{k}}, \varphi (t+p_1(z_0))
\right\rangle
\int\limits_{\mathbb{S}^{2n-1}}\Omega_{0}^{+}(\theta)|\mathfrak{p}_{k}(\theta
)|^{-\frac{2n}{k}}d\theta,
\end{equation}
for the set of directions where $\mathfrak{p}_k(\theta)>0$ and by
:
\begin{equation} \label{negative}
I_{-}(\lambda) \sim \frac{1}{k} \lambda^{-\frac{2n}{k}}
\left\langle |t|_{-}^{\frac{2n-k}{k}}, \varphi (t+p_1(z_0))
\right\rangle
\int\limits_{\mathbb{S}^{2n-1}}\Omega_{0}^{-}(\theta)|\mathfrak{p}_{k}(\theta
)|^{-\frac{2n}{k}}d\theta,
\end{equation}
for the directions where $\mathfrak{p}_k (\theta)$ is negative. \medskip\\
\newline\textbf{Microlocal contribution of the set $\mathfrak{p}_k(\theta)=0$.}\newline%
\textit{Case of $k>2n$.}\\
Here we examine the contribution of terms :
\begin{equation*}
\int e^{\frac{i}{h}\chi _{0}\chi _{1}^{k}\chi_{2}}A_{i}(\chi
_{0},\chi _{1},\chi_{2})d\chi _{0}d\chi _{1}d\chi_{2},
\end{equation*}
but Lemma \ref{Theo IO 2eme carte} shows that these are given by :
\begin{equation*}
\frac{1}{k}\Gamma(\frac{2n}{k})h^{\frac{2n}{k}}\int |\chi_0 \chi_2
|^{-\frac{2n}{k}} \exp(i\frac{\pi n}{k}\mathrm{sign}(\chi_0
\chi_2) )\tilde{A}_{i}(\chi_0,0,\chi_2)d\chi_0
d\chi_2+\mathcal{O}(h ^{\frac{2n+1}{k}}\mathrm{log}^2(h)).
\end{equation*}
Writing the delta-Dirac distribution as an oscillatory integral
leads to :
\begin{equation*}
\frac{1}{k}\Gamma(\frac{2n}{k})h^{\frac{2n}{k}} \frac{1}{2\pi}\int
e^{iz\chi_1}|\chi_0 \chi_2 |^{-\frac{2n}{k}}\exp(i\frac{\pi
n}{k}\mathrm{sign}(\chi_0 \chi_2) )
\tilde{A}_{i}(\chi_0,\chi_1,\chi_2)d\chi_0 d\chi_1d\chi_2dz,
\end{equation*}
with the amplitude :
\begin{equation*}
\tilde{A}_{i}(\chi_0,\chi_1,\chi_2)=\int
\chi^{\star}(\Omega_i(\theta) a(t,r\theta)
|J\chi|)d\chi_3...d\chi_{2n}.
\end{equation*}
Now, we return to the initial coordinates and with $\chi_1 =r$
this gives :
\begin{equation*}
\frac{1}{k}\Gamma(\frac{2n}{k})h^{\frac{2n}{k}} \frac{1}{2\pi}\int
e^{izr}|(\chi_0 \chi_2)(t,r,\theta)|^{-\frac{2n}{k}}
\Omega_i(\theta) \exp(i\frac{\pi n}{k}\mathrm{sign}(\chi_0 \chi_2)
)a(t,r\theta)dtdrd\theta dz,
\end{equation*}
if we use that
$(\chi_0,\chi_1,\chi_2)(t,0,\theta)=(t,0,\mathfrak{p}_k(\theta))$,
we obtain :
\begin{equation*}
\frac{1}{k}\Gamma(\frac{2n}{k})h^{\frac{2n}{k}}\int\limits_{\mathbb{R}\times\mathbb{S}^{2n-1}}
|t \mathfrak{p}_{k}(\theta) |^{-\frac{2n}{k}}\exp(i\frac{\pi
n}{k}\mathrm{sign}(t\mathfrak{p}_k(\theta)))
\Omega_i(\theta)a(t,0)dt d\theta.
\end{equation*}
Hence, when $\mathfrak{p}_k(\theta)$ is positive, we have :
\begin{gather*}
\frac{1}{k}\Gamma(\frac{2n}{k})\int\limits_{\mathbb{R}}
|t|^{-\frac{2n}{k}} \exp(i\frac{\pi n}{k}\mathrm{sign}(t))
\hat{\varphi}(t) e^{itp_1(z_0)}dt \\
= \frac{1}{k}\Gamma(\frac{2n}{k})\left\langle
|t|_{+}^{-\frac{2n}{k}} \exp(i\frac{\pi
n}{k})+|t|_{-}^{-\frac{2n}{k}} \exp(-i\frac{\pi
n}{k}),\hat{\varphi}(t) e^{itp_1(z_0)} \right\rangle .
\end{gather*}
Similarly, when $\mathfrak{p}_k(\theta)$ is negative we have :
\begin{gather*}
\frac{1}{k}\Gamma(\frac{2n}{k})\int\limits_{\mathbb{R}}
|t|^{-\frac{2n}{k}} \exp(-i\frac{\pi n}{k}\mathrm{sign}(t))
\hat{\varphi}(t) e^{itp_1(z_0)}dt \\
= \frac{1}{k}\Gamma(\frac{2n}{k})\left\langle
|t|_{+}^{-\frac{2n}{k}} \exp(-i\frac{\pi
n}{k})+|t|_{-}^{-\frac{2n}{k}} \exp(i\frac{\pi
n}{k}),\hat{\varphi}(t) e^{itp_1(z_0)} \right\rangle .
\end{gather*}
Hence, by summation the contribution is :
\begin{gather*}
\frac{1}{k}\Gamma(\frac{2n}{k})h^{\frac{2n}{k}}\left(\cos(\frac{\pi
n}{k}) \left\langle |t|^{-\frac{2n}{k}},\hat{\varphi}(t)
e^{itp_1(z_0)}\right\rangle \int\limits_{\mathbb{S}^{2n-1}}
\Omega_i (\theta) |\mathfrak{p}_k (\theta)
|^{-\frac{2n}{k}}d\theta \right.\\
\left. +i\sin(\frac{\pi n}{k})\left\langle
|t|^{-\frac{2n}{k}}\mathrm{sign}(t),\hat{\varphi}(t)
e^{itp_1(z_0)}\right\rangle \int\limits_{\mathbb{S}^{2n-1}}
\Omega_i (\theta)|\mathfrak{p}_k (\theta)
|^{-\frac{2n}{k}}\mathrm{sign}(\mathfrak{p}_k
(\theta))d\theta\right ).
\end{gather*}
If we use the classical relations :
\begin{gather*}
\mathcal{F}(|x|^\lambda )(\xi)=-2\sin(\frac{\lambda \pi}{2}) \Gamma(\lambda +1) |\xi|^{-\lambda-1},\\
\mathcal{F}(|x|^\lambda \mathrm{sign}(x))(\xi)=2i
\cos(\frac{\lambda \pi}{2}) \Gamma(\lambda +1) |\xi|^{-\lambda-1}
\mathrm{sign}(\xi),
\end{gather*}
we obtain, after some manipulations, that the contribution is :
\begin{gather*}
\frac{1}{k}h^{\frac{2n}{k}}\left( \left\langle
|t|^{\frac{2n}{k}-1},\varphi(t+p_1(z_0))\right\rangle
\int\limits_{\mathbb{S}^{2n-1}} \Omega_i (\theta) |\mathfrak{p}_k
(\theta) |^{-\frac{2n}{k}}d\theta\right.\\
\left. + \left\langle
|t|^{\frac{2n}{k}-1}\mathrm{sign}(t),\varphi(t+p_1(z_0))\right\rangle
\int\limits_{\mathbb{S}^{2n-1}} \Omega_i (\theta)|\mathfrak{p}_k
(\theta) |^{-\frac{2n}{k}}\mathrm{sign}(\mathfrak{p}_k
(\theta))d\theta \right).
\end{gather*}
We can split the integral with respect to $d\theta$ into two parts
to finally obtain :
\begin{gather*}
\frac{1}{k}h^{\frac{2n}{k}}\left( \left\langle
|t|_{+}^{\frac{2n}{k}-1},\varphi(t+p_1(z_0))\right\rangle
\int\limits_{\mathbb{S}^{2n-1}\cap \{\mathfrak{p}_k \geq 0\}}
\Omega_i (\theta) |\mathfrak{p}_k
(\theta) |^{-\frac{2n}{k}}d\theta\right.\\
\left. + \left\langle
|t|_{-}^{\frac{2n}{k}-1},\varphi(t+p_1(z_0))\right\rangle
\int\limits_{\mathbb{S}^{2n-1}\cap \{\mathfrak{p}_k \leq 0\}}
\Omega_i (\theta)|\mathfrak{p}_k (\theta) |^{-\frac{2n}{k}}d\theta
\right).
\end{gather*}
With Eq. (\ref{positive}) and (\ref{negative}), by summation on
the partition of unity the main contribution to the trace formula
is given by :
\begin{gather}
\gamma _{z_{0}}(E_{c},h)\simeq\frac{1}{k}\frac{h^{\frac{2n}{k}-n}}{(2\pi)^{n}}%
\left ( \left\langle
|t|_{+}^{\frac{2n}{k}-1},\varphi(t+p_1(z_0))\right\rangle
\int\limits_{\mathbb{S}^{2n-1}\cap \{ \mathfrak{p}_k \geq 0\}}
|\mathfrak{p}_k
(\theta) |^{-\frac{2n}{k}}d\theta\right. \notag \\
\left. + \left\langle
|t|_{-}^{\frac{2n}{k}-1},\varphi(t+p_1(z_0))\right\rangle
\int\limits_{\mathbb{S}^{2n-1}\cap \{ \mathfrak{p}_k \leq 0\}}
|\mathfrak{p}_k (\theta) |^{-\frac{2n}{k}}d\theta \right ).
\label{final}
\end{gather}
And this proves the first statement of Theorem \ref{Main1} for
$k>2n$. $\hfill{\blacksquare}$\medskip\\
\textit{Case of $k=2n$.}\\
Here the contribution is given by :
\begin{equation*}
\int\limits_{\mathbb{R\times R}_{+}\times
\mathbb{S}^{2n-1}}e^{\frac{i}{h}\Psi (t,r,\theta
)}\Omega_{i}(\theta) a(t,r\theta )r^{2n-1}dtdrd\theta =\int
e^{\frac{i}{h}\chi _{0}\chi _{1}^{k}\chi_{2}}A_{i}(\chi _{0},\chi
_{1},\chi_{2})d\chi _{0}d\chi _{1}d\chi_{2},
\end{equation*}
and we recall that the amplitude is given by :
\begin{equation*}
A_{i}(\chi _{0},\chi _{1},\chi_{2})=\int \chi ^{\ast
}(\Omega_{i}(\theta)a(t,r\theta )r^{2n-1}|J\chi |)d\chi
_{3}...d\chi _{2n}.
\end{equation*}
With $A_{i}=\chi_{1}^{2n-1}\tilde{A}_{i}$, from Lemma \ref{Theo IO
2eme with log} we know that we have :
\begin{equation*}
\int e^{\frac{i}{h}\chi _{0}\chi
_{1}^{k}\chi_{2}}\tilde{A}_{i}(\chi _{0},\chi
_{1},\chi_{2})\chi_{1}^{2n-1}d\chi _{0}d\chi
_{1}d\chi_{2}=\frac{\pi h}{n} \log (h)
\tilde{A}_{i}(0,0,0)+\mathcal{O}(h).
\end{equation*}
Where, by construction, we have defined :
\begin{equation*}
\tilde{A}_{i}(\chi _{0},\chi _{1},\chi_{2})=\chi_1^{-(2n-1)}\int
\chi ^{\ast }(\Omega_{i}(\theta)a(t,r\theta )r^{2n-1}|J\chi
|)d\chi _{3}...d\chi _{2n}.
\end{equation*}
With $z=(z_1,z_2,z_3)$, we write the amplitude as :
\begin{equation*}
\tilde{A}_{i}(0,0,0)=\frac{1}{(2\pi)^3}\int\limits
e^{-i\left\langle z, (\chi_0,\chi_1,\chi_2)\right\rangle}
\tilde{A}_{i}(\chi_0,\chi_1,\chi_2)dzd\chi_0 d\chi_1 d\chi_2,
\end{equation*}
and by inversion of our diffeomorphism we obtain :
\begin{gather*}
\tilde{A}_{i}(0,0,0)=\frac{1}{(2\pi)^3}\int\limits
e^{-i\left\langle z, (\chi_0,\chi_1,\chi_2)\right\rangle}
\Omega_{i}(\theta)a(t,r\theta )drdtd\theta dz.
\end{gather*}
Using that $(\chi_0,\chi_1)=(t,r)$, by integration w.r.t.
$(t,r,z_1,z_2)$ we get :
\begin{equation*}
\tilde{A}_{i}(0,0,0)=\frac{1}{2\pi}\int\limits e^{-i
z_3\chi_2(0,0,\theta)} \Omega_{i}(\theta)a(0,0)d\theta dz_3.
\end{equation*}
By construction $\chi_2(t,0,\theta)=\mathfrak{p}_k(\theta)$ and,
since $\mathfrak{p}_k(\theta)$ is an admissible coordinate on
$\mathrm{supp}(\Omega_{i})$, we can use the change of coordinates
$u=\mathfrak{p}_k(\theta)$, this leads to :
\begin{equation*}
\tilde{A}_{i}(0,0,0)=\frac{1}{2\pi}a(0,0)\int\limits e^{-i z_3 u}
\Omega_{i}(\theta)dL_{\mathfrak{p}_k}(\theta) dudz_3,
\end{equation*}
with $d\mathfrak{p}_k \wedge dL_{\mathfrak{p}_k}(\theta)=d\theta$.
If we introduce :
\begin{equation*}
C_{\mathfrak{p}_k}=\{\theta \in \mathbb{S}^{2n-1} \text{ / }
\mathfrak{p}_k(\theta)=0 \},
\end{equation*}
the sum over the partition of unity gives :
\begin{equation*}
\sum\limits_{i}\int\limits_{C_{\mathfrak{p}_k}} \Omega_i (\theta)
dL_{\mathfrak{p}_k}(\theta)=\mathrm{LVol}(C_{\mathfrak{p}_k}),
\end{equation*}
where $\mathrm{LVol}$ is the Liouville volume attached to
$\mathfrak{p}_k(\theta)$ on $\mathbb{S}^{2n-1}$. Finally, since we
have $a(t,0)=\hat{\varphi}(t)\mathrm{exp}(itp_1(z_0))$, the second
statement of Theorem \ref{Main1} holds.

\end{document}